\documentclass[a4paper,12pt]{amsart}
\usepackage{amsmath,inputenc,amssymb}
\title[Blow-up for NLS on domains]{Remarks
on the blow-up for the Schr{\"o}dinger equation with critical
mass on a plane domain}
\author{Valeria Banica}

\newtheorem{lemma}{Lemma}[section]
\newtheorem{theorem}[lemma]{Theorem}
\newtheorem{prop}[lemma]{Proposition}
\newtheorem{corollary}[lemma]{Corollary}

\newtheorem{remark}[lemma]{Remark}
\newtheorem{definition}[lemma]{Definition}

\newcommand{\nL}[1]{\|#1 \|_2}
\newcommand{\mo}[1]{|#1|}
\newcommand{\ls}{\lesssim}

\newcommand{\al}{\alpha}

\newcommand{\la}{\lambda}
\newcommand{\DP}[1]{{\partial_ #1}}

\def\tend{{\rightarrow}}
\newcommand{\R}{\mathbb{R}}

\begin{document}

\begin{abstract}
In this paper we concentrate on the analysis of the critical mass 
blowing-up solutions for the cubic focusing Schr{\"o}dinger 
equation with Dirichlet boundary conditions, posed on a plane domain. 
We bound the blow-up rate from below, for bounded and unbounded domains. 
If the blow-up occurs on the boundary, the blow-up rate is proved to grow 
faster than $(T-t)^{-1}$, the expected one. Moreover, we show that 
blow-up cannot occur on the boundary, under certain geometric conditions 
on the domain.
\end{abstract}

\address{V. Banica, Universita di Pisa, Dipartimento di Matematica, 
Via Buonarroti 2, 56127 Pisa, Italy}
\email{Valeria.Banica@mail.dm.unipi.it}

\maketitle

\section{Introduction} 

Let us first recall the known results for the $\mathbb{R}^n$ case. \\
Consider the nonlinear Schr{\"o}dinger equation on $\mathbb{R}^n$, for 
$p\geq 1$,
$$(S)\left\{\begin{array}{c} 
i\DP{t}u+\Delta u+\mo{u}^{p-1}u=0,\\ 
u(0)=u_0.
\end{array}\right.$$ 
The associated Cauchy problem is locally well posed in $\mathbb{H}^1$ for 
$p<1+\frac{4}{n-2}$ (\cite{GV23}, \cite{Ka3}).

The Gagliardo-Nirenberg inequality
$$\|v\|_{p+1}^{p+1}\leq C_{p+1} \|v\|_2^{2+(p-1)\frac{2-n}{2}}
\|\nabla v\|_2^{(p-1)\frac{n}{2}}$$
implies that the energy of the solution $u$ of the equation $(S)$,
$$E(u)=
\frac{1}{2}\int_{\mathbb{R}^n}\mo{\nabla u}^2dx-
\frac{1}{p+1}\int_{\mathbb{R}^n}\mo{u}^{p+1}dx,$$
is bounded from below by
$$\|\nabla u\|_2^2\left(\frac{1}{2}-
\frac{C_{p+1}}{p+1}\|u\|_2^{2+(p-1)\frac{2-n}{2}}
\|\nabla u\|_2^{(p-1)\frac{n}{2}-2}\right).$$
As a consequence, if $p<1+\frac{4}{n}$, since the mass is conserved, 
the gradient of $u$ is controlled by the energy. 
Therefore the solution does not blow up and global existence occurs.

The power $p=1+\frac{4}{n}$ is a critical power, in the sense that the 
nonlinearity is strong enough to generate solutions blowing up in a 
finite time.
However, even in this case, we have a global result for small initial 
conditions.\\
Indeed, in the case $p=1+\frac{4}{n}$, if the mass of the initial 
condition is small enough so that
$$\frac{C_{2+\frac{4}{n}}}{2+\frac{4}{n}}\|u\|_2^{\frac{4}{n}}<
\frac{1}{2},$$
then the energy controls the gradient and again, the global existence is 
proved for the equation $(S)$.\\
For this particular value of $p$, Weinstein has given a sharpening of the
Gagliardo-Nirenberg inequality (\cite{We3}).
By variational methods using Lions concentration-compacity lemma 
(\cite{Li13}, \cite{Li23}), he 
obtained the existence of a minimizer $Q$ for the optimal constant of 
Gagliardo-Nirenberg's inequality
$$\frac{1}{C_{2+\frac{4}{n}}}=
\inf_{v\in\mathbb{H}^1(\mathbb{R}^n)}
\frac{\|v\|_2^{\frac{4}{n}}
\|\nabla v\|_2^2}{\|v\|_{2+\frac{4}{n}}^{2+\frac{4}{n}}}.$$
This minimizer satisfies the Euler-Lagrange equation
$$\Delta Q+Q^{1+\frac{4}{n}}=\frac{2}{n}\,Q.$$
Such a positive function, called ground state of the nonlinear 
Schr{\"o}dinger equation, is radial, exponentially decreasing at 
infinity and regular.
Recently, Kwong has shown that it is unique up to a translation 
(\cite{Kv3}).
Moreover, it verifies Pohozaev's identities
$$\left\{\begin{array}{c}
\|\nabla Q\|_2^2-\|Q\|_{2+\frac{4}{n}}^{2+\frac{4}{n}}+
\frac{2}{n}\|Q\|_2^2=0,\\
(n-2)\|\nabla Q\|_2^2-\frac{n^2}{n+2}\|Q\|_{2+
\frac{4}{n}}^{2+\frac{4}{n}}+2\|Q\|_2^2=0,
\end{array}\right.$$
which lead to the following relations between the norms of $Q$
\begin{equation}\label{normQ}
\left\{\begin{array}{c}
\|Q\|_{2+\frac{4}{n}}^{2+\frac{4}{n}}=\frac{n+2}{n}\|Q\|_2^2,\\
\|\nabla Q\|_2^2=\|Q\|_2^2.
\end{array}\right.
\end{equation}
Then the optimal value for the constant of the Gagliardo-Nirenberg 
inequality is 
$$C_{2+\frac{4}{n}}=\frac{n+2}{n}\frac{1}{\|Q\|_2^{\frac{4}{n}}}.$$
In conclusion, if $p=1+\frac{4}{n}$, the solutions of the equation 
$(S)$ with initial condition of mass smaller than the one of the ground 
state 
$$\|u\|_2<\|Q\|_2,$$
are global in time.

The mass $\|Q\|_2$ is critical, in the sense that we can construct
as follows solutions of mass equal to $\|Q\|_2$, which blows up in 
finite time. 
Since $p=1+\frac{4}{n}$, the pseudo-conformal transform of a 
solution $u$ of $(S)$
$$\frac{1}{t^{\frac{n}{2}}}e^{i\frac{\mo{x}^2}{4t}}
u\left(-\frac{1}{t},\frac{x}{t}\right),$$
is also a solution of $(S)$ (\cite{Ca03}). 
So, from a stationary solution on $\mathbb{R}^n$
$$e^{it}Q(x),$$
for all positive $T$,
$$u(t,x)=\frac{e^{\frac{i}{T-t}}}{(T-t)^{\frac{n}{2}}}\,
e^{-i\frac{|x|^2}{4(T-t)}}\,Q\left(\frac{x}{T-t}\right),$$
is a solution blowing up at the time $T$. 
Moreover, Merle proved in \cite{Me3} that all blowing-up solutions 
on $\mathbb{R}^n$ with critical mass $\nL{Q}$ are of this type, up to the 
invariants of the equation. 
The proof is based on a result of concentration of Weinstein 
(\cite{We23}, see Lemma \ref{Wein}) and on the study of the first
order momentum
$$f(t)=\int_{\mathbb{R}^n}|u(t,x)|^2xdx,$$
and of the virial
$$g(t)=\int_{\mathbb{R}^n}|u(t,x)|^2|x|^2dx$$
associated to a solution $u$ of the equation $(S)$. 
The conservative properties of these two quantities on $\mathbb{R}^n$, 
in the case of the critical power $1+\frac{4}{n}$, play an important 
role in Merle's proof. 
The derivative of the first order momentum is constant in time
$$\partial_t^2 f=0,$$
and $g$ satisfies the virial identity (\cite{Ca03})
$$\partial_t^2 g=
16E(u)-4\frac{n(p-1)-4}{p+1}\int_{\mathbb{R}^n}|u|^{p+1}dx=16E(u).$$

In certain cases of initial conditions with mass larger than $\|Q\|_2$ 
recent achievements were done by Merle and Rapha{\"e}l, concerning the
blow-up rate and the blow-up profile (\cite{MR13}, \cite{MR23}).

For the equation $(S_p)$ with $p\geq 1+\frac{4}{n}$, Zakharov \cite{Za3} 
and Glassey \cite{Gl3} had obtained that the solutions of negative 
energy are blowing up in finite time. 
The same result for solutions of nonnegative 
energy is valid under certain conditions on the derivatives of the 
virial (\cite{SS3}). 
The proof is based on an upper bound of the virial in terms of its 
first and second derivative, which implies the cancellation of the 
virial at a finite time $T$. 
Since the mass is conserved, it follows that the solution must blow up 
at the time $T$.

\bigskip

In this paper we are concerned with the nonlinear Schr{\"o}dinger 
equation posed on a regular domain $\Omega$ of $\R^n$, with Dirichlet 
boundary conditions
$$\left\{\begin{array}{c} 
i\DP{t}u+\Delta u+\mo{u}^{p-1}u=0,\\
u_{|\R\times\partial\Omega}=0,\\
u(0)=u_0.
\end{array}\right.$$
The Cauchy problem is locally well posed on 
$\mathbb{H}^2\cap\mathbb{H}^1_0(\Omega)$ in dimension $2$ and $3$. 
In dimension $2$, for nonlinearities less than cubic, 
Vladimirov \cite{Vl3} and Ogawa and Ozawa \cite{OO3} have shown 
the well-posedness of the Cauchy problem on 
$\mathbb{H}^1_0(\Omega)$, but without the uniform continuity of the 
flow on bounded sets of $\mathbb{H}^1_0(\Omega)$. 
For nonlinearities stronger than cubic in dimension $2$, 
or for any power nonlinearity $p$, in dimension higher than $2$, 
the Cauchy problem on $\mathbb{H}^1_0(\Omega)$ is open.

For the equation with power $p<1+\frac{4}{n}$, one can show as for the 
case $\R^n$ that the $\mathbb{H}^1_0(\Omega)$ solutions are global in 
time.
For the equation with power $p\geq 1+\frac{4}{n}$, posed on a star-shaped 
domain of $\R^n$, Kavian has proved the blow-up in finite time of the 
$\mathbb{H}^2\cap\mathbb{H}^1_0(\Omega)$ 
solutions of negative energy or of positive energy but under some 
conditions on the first and second derivatives of the virial 
(\cite{Kav3}).
His proof follows the one on $\R^n$ (\cite{Gl3}), by 
estimating via the geometric condition on $\Omega$ the boundary terms 
which appear in the second derivative of the virial.

From now on we will analyze the cubic equation on $\Omega$
$$(S_\Omega)\,\,\,\,\,\left\{\begin{array}{c} 
i\DP{t}u+\Delta u+\mo{u}^2u=0,\\
u_{|\R\times\partial\Omega}=0,\\
u(0)=u_0.
\end{array}\right.$$
Let us first notice that the conservations of the mass and of the
energy of the solutions are still valid. 
The Cauchy problem is locally well posed on 
$\mathbb{H}^2\cap\mathbb{H}^1_0(\Omega)$, and also on 
$\mathbb{H}^1_0(\Omega)$ apart from the property of uniform continuity of 
the flow, not known to hold. 
The usual Strichartz inequalities are no longer valid and the loss of 
derivatives is stronger than in the case of a compact manifold 
(\cite{BGT33}).

As in the case of the plane, for initial conditions with mass 
smaller than the one of the ground state, the Cauchy problem 
is globally well-posed on $\mathbb{H}^2\cap\mathbb{H}^1_0(\Omega)$. 
The proof, given by Br{\'e}zis and Gallou{\"e}t, is based on logarithmic 
type estimates (\cite{BG3}). 
This result has been extended to the natural space 
$\mathbb{H}^1_0(\Omega)$, apart from the 
uniform continuity of the flow (\cite{Vl3},\cite{OO3},\cite{Ca03}).

The critical mass for blow-up is $\|Q\|_2$, as in the case of the 
equation posed on $\R^2$. 
More precisely, the following result holds.
\begin{theorem}\label{ex} (Burq-G{\'e}rard-Tzvetkov \cite{BGT33})
Let $\Omega$ be a regular bounded domain of $\R^2$. 
Let $x_0\in\Omega$ and $\psi\in\mathcal{C}^\infty_0(\Omega)$ be a function 
equal to $1$ near $x_0$.
Then there exist positive numbers $\kappa$ and $\alpha_0$ such 
that for all $\alpha>\alpha_0$, there exists a time $T_\alpha$ and a 
function $r_\alpha$ defined on $[0,T_\alpha[\times\Omega$, satisfying 
$$\|r_\alpha(t)\|_{\mathbb{H}^2(\Omega)}\leq 
ce^{-\frac{\kappa}{T_\alpha-t}},$$
such that 
$$u(t,x)=\psi(x)\frac{e^{\frac{i}{\alpha^2(T_\alpha-t)}}}
{\alpha(T_\alpha-t)}e^{-i\frac{|x-x_0|^2}{4\alpha(T_\alpha-t)}}
Q\left(\frac{x-x_0}{\alpha(T_\alpha-t)}\right)+r_\alpha(t,x),$$
is a critical mass solution of $(S_\Omega)$, blowing up at $x_0$ at 
the time $T_\alpha$ with the blow-up rate $\frac{1}{T_\alpha-t}$.
\end{theorem}

The proof, following an idea of Ogawa and Tsutsumi (\cite{OT3}), 
is based on a fixed point method which allows to complete 
the cut-off of the explicit blowing up solution on $\R^2$ at $x_0$ 
to a blowing up solution on $\Omega$ at $x_0$. 
Theorem 4.1.1 implies in particular that at every point of $\Omega$ there are 
explosive solutions. 
Moreover, the proof is still valid for the torus $\mathbb{T}^2$ 
and for a larger class of subsets of the plane, which satisfy the
property of 2-continuation, from
$\mathbb{H}^2\cap\mathbb{H}_0^1(\Omega)$ to $\mathbb{H}^2(\R^2)$, and 
for which the Laplacian domain
$$D(-\Delta_{\Omega})=
\{u\in\mathbb{H}_0^1(\Omega), \Delta u\in\mathbb{L}^2(\Omega)\},$$
is $\mathbb{H}^2\cap\mathbb{H}^1_0$. 
Such subsets are for example the domains with compact regular boundary 
and convex polygons bounded or unbounded.

As in the $\R^n$ case, the following lemma, due to Weinstein, 
will give us the general behavior of a blowing-up solution of critical
mass on a domain.
\begin{lemma}\label{Wein} (Weinstein \cite{We23})
Let $u_k\in\mathbb{H}^1(\mathbb{R}^n)$ be a sequence of functions of 
critical mass satisfying
$$\left\{\begin{array}{c}
\beta_k=\|\nabla u_k\|_2\underset{k\tend \infty}{\longrightarrow}\infty,\\
E(u_k)\underset{k\tend \infty}{\longrightarrow}c<\infty.
\end{array}\right.$$
Then there exist points $x_k\in\R^d$ and $\theta_k\in\R$ such that 
in $\mathbb{H}^1(\R^n)$
$$\frac{e^{i\theta_k}}{\beta_k^{\frac{n}{2}}}
u_k\left(\frac{x}{\beta_k}+x_n\right)
\underset{k\tend \infty}{\longrightarrow}
\frac{1}{\omega^{\frac{n}{2}}}Q\left(\frac{x}{\omega}\right),$$
where $\omega=\|\nabla Q\|_2$.
\end{lemma}

Let $u$ be a solution of $(S_\Omega)$ that blows up at the finite 
time $T$, that is 
$$\la(t)=
\frac{\nL{\nabla u(t)}}{\nL{\nabla Q}}
\underset{t\tend T}{\longrightarrow}\infty.$$
Consider $u$ to be extended by zero outside $\Omega$. 
By combining Lemma \ref{Wein} for families $u_k=u(t_k)$ with $t_k$ 
sequences 
convergent to $T$ with the result of Kwong on the uniqueness of the 
ground state (\cite{Kv3}), there exist $\theta(t)$ real numbers and 
$x(t)\in\R^2$ such that in $\mathbb{H}^1(\mathbb{R}^2)$ 
\begin{equation}\label{limW}
\frac{e^{i\theta(t)}}{\la(t)}
u\left(t,\frac{x}{\la(t)}+x(t)\right)
\underset{t\tend T}{\longrightarrow}Q(x).
\end{equation}
Then, in the space of distributions, 
\begin{equation}\label{limx}
\mo{u(t,\cdot+x(t))}^2
\underset{t\tend T}{\longrightarrow} \nL{Q}^2\,\delta_0.
\end{equation}

\bigskip

In this paper we concentrate on the further analysis of the blowing-up 
solutions with critical mass on a plane domain. 
The results are the following. 

\begin{theorem}\label{art31} 
Let $u$ be a $\mathcal{C}([0,T[,\mathbb{H}_0^1)$ solution of the 
Schr{\"o}dinger equation $(S_\Omega)$,
which has critical mass and blows up at the finite time $T$. 

i) For bounded domains, the blowing-up rate is lower bounded by 
$$\frac{1}{T-t}\lesssim\|\nabla u(t)\|_2.$$

ii) If there exist solutions $u$ of critical mass blowing up at a 
finite time $T$ on the boundary of $\Omega$, that is if the 
concentration parameter $x(t)$ converges as $t\tend T$ to a point on the 
boundary, then the blowing-up rate satisfies
$$\underset{t\tend T}{\lim}(T-t)\|\nabla u(t)\|_2=\infty.$$
\end{theorem}

The main difficulty for the Schr{\"o}dinger equation posed on a
domain is that the conservation of the derivative of the first 
momentum and the virial identity fail.\par 
In order to avoid this difficulty, we shall use systematically in the 
proof of Theorem \ref{art31} a Cauchy-Schwarz type inequality derived 
from Weinstein's inequality. 
Precisely, we show that if $v$ is a $\mathbb{H}^1(\mathbb{R}^2)$ function 
of critical or subcritical mass, then
$$\left\vert\int\Im{(v\nabla\overline{v}\,)}\nabla\theta dx\right\vert
\leq 
\left(2E(v)\int\mo{v}^2\mo{\nabla\theta}^2dx\right)^{\frac{1}{2}}$$ 
for all real function $\theta$. 
This inequality allows us to estimate the virial, that we shall assume 
to be localized if $\Omega$ is unbounded (see Remark \ref{unbound}). 
The lower bound for the blowing-up rate is the same as the one found by 
Antonini on the torus (\cite{An3}).

By following the approach of Weinstein in \cite{We33}, and the recent 
results of Maris in \cite{M3}, we analyze the convergence
to the ground state of the modulations of the solutions \eqref{limW},
and we obtain, for bounded domains, the following additional informations.

\begin{prop}\label{equivir}
i) The blow-up rate verifies
$$\int_{\Omega} |u(t,x)|^2|x-x(t)|^2dx\approx\frac{1}{\|\nabla u(t)\|_2^2}.$$
ii) The concentration parameter $x(t)$ can be chosen to be as 
the first order momentum 
$$x(t)=\frac{\int_{\Omega} |u(t,x)|^2xdx}{\|Q\|_2^2}.$$
\end{prop}

\begin{corollary}\label{radial}
If $\Omega$ is a disc centered at $0$ and if equation $(S_\Omega)$ is 
considered to be invariant under rotations, then $x(t)$ can be chosen to be 
$0$, and we have
$$g(t)\approx\frac{1}{\|\nabla u(t)\|_2^2}.$$
\end{corollary}
\begin{remark}\label{unbound}
For unbounded domains, if the
solution concentrates at one point, that is if $x(t)$ converges as 
$t\tend T$, then the first assertion of Theorem \ref{art31} is true, and 
so are the assertions of Proposition \ref{equivir},
for the virial and the first order momentum localized at the blow-up 
point (see \S4).
\end{remark}

There is no known example of a solution of nonlinear Schr{\"o}dinger 
equation with a blow-up rate larger than $\frac{1}{T-t}$, neither in the 
case of supercritical mass, nor in the case of supercritical 
nonlinearities.

Therefore we expect that the blowing-up rate grows exactly like 
$\frac{1}{T-t}$ 
and that the profiles are the ones on $\mathbb{R}^2$ modulo an 
exponentially decreasing in $\mathbb{H}^1$ function. 

Since it is not likely that the blowing-up rate at the boundary grows 
strictly faster than $\frac{1}{T-t}$, we also expect that there 
are no solutions blowing-up on the boundary of a domain. 
This is confirmed for certain simple cases by the following result.
\begin{theorem}\label{art32}
If $\Omega$ is a half-plane or a plane sector, then there are 
no solutions blowing-up in a finite time on the boundary of the 
half-plane or in the corner of the sector respectively.
\end{theorem}
 
Indeed, under these geometric hypotheses on $\Omega$, the boundary terms 
which appear in the second derivative of the virial associated to a 
blowing-up solution of $(S_\Omega)$ cancel, so we have, as on $\R^n$, 
the virial identity 
$$\partial^2_t g=16E(u).$$
The proof then follows the one by Merle in \cite{Me3} for the equation 
posed on $\R^n$, and we obtain that all explosive solutions on 
$\Omega$ must be of the type
$$\frac{e^{\frac{i}{T-t}}}{(T-t)^{\frac{n}{2}}}\,
e^{-i\frac{|x|^2}{4(T-t)}}\,Q\left(\frac{x}{T-t}\right),$$
up to the invariants of the equation.
Therefore we arrive at a contradiction by looking at the support of the 
solution. 

\begin{remark}The results of this paper are valid also in higher dimension, 
but are presented in the 2-dimensional case. 
The reason is that the only existence result of blowing-up solutions with 
critical mass on a domain is Theorem \ref{ex}.
\end{remark}

\begin{remark}
In the case of Neumann boundary conditions, the situation changes radically
for blow up at the boundary. 
For blow up inside the domain, the situation is expected to be similar to 
the Dirichlet case.\par
Blow-up can appear at the boundary of a half-plane, as it can easily be seen 
by taking the restriction to the half-plane of the explicit solution on 
$\mathbb{R}^2$, blowing-up at zero. 
Notice also that this is a blowing-up solution of mass $\frac{\|Q\|_2^2}{2}$. 
By the same method one can construct solution of mass 
$\frac{\|Q\|_2^2}{n}$, blowing-up at the corner of a plane sector of angle 
$\frac{2\pi}{n}$. 
So for unbounded domains the critical mass for blow up on the boundary 
probably depends on the domain.

By using the techniques of Theorem \ref{ex}, blowing-up solution of this mass 
$\frac{\|Q\|_2^2}{2}$ can be contructed to blow up on the boundary of a 
bounded set, in a point around which the boundary is locally a line. 
For bounded domains, we expect that the critical mass for blow up on the 
boundary should be $\frac{\|Q\|_2^2}{2}$.

One of the difficulties is that the best constant in the corresponding 
Gagliardo-Nirenberg inequality 
$$\|v\|_4^4\leq A\|v\|_2^2(\|\nabla v\|_2^2+B\|v\|_2^2),$$
is not known for the case of Neumann boundary conditions.
\end{remark}

The paper is organized as follows. 
Section \S2 contains some results on general domains. 
We prove a Cauchy-Schwarz type inequality for critical and subcritical 
mass functions, which we will use to show Theorem \ref{art31}. 
The nature of the convergence to the ground state of the modulations of
the solutions is analyzed, by spectral theory techniques given in the 
Appendix. These concentration results will be
used later to prove Theorem \ref{art32} and Proposition \ref{equivir}.  
Moreover, we calculate the derivatives in time for a virial type function.
In \S3, by studying the virial, the lower-bound of the blowing-up rate 
is proved for bounded domains $\Omega$. 
In this section, we also give the proof of Proposition \ref{equivir}.
In \S4, by introducing a localized virial, we find the same lower-bound
for the blowing-up rate for unbounded domains. 
Section \S5 contains the results regarding the explosion on 
the boundary of $\Omega$.

{\bf{Acknowledgment.}} 
I would like to thank my advisor Patrick G{\'e}rard for having introduced me 
to this beautiful subject and for having guided this work. 
I am also grateful to Mihai Maris for useful discussions concerning the 
Appendix.
 
\section{Results on general domains}
\subsection{A Cauchy-Schwarz inequality for subcritical mass 
functions}\label{CS} 
\begin{lemma} 
Let $\theta$ be a real valued function. 
All $v\in\mathbb{H}^1(\mathbb{R}^2)$ with critical or subcritical 
mass satisfy 
$$(*)\left\vert\int\Im{(v\nabla\overline{v}\,)}\nabla\theta dx
\right\vert\leq 
\left(2E(v)\int\mo{v}^2\mo{\nabla\theta}^2dx\right)^{\frac{1}{2}}.$$ 
\end{lemma} 

\begin{proof} 
The precised version of the Gagliardo-Niremberg inequality, presented 
in the introduction, is, for function $w$ in $\mathbb{H}^1(\mathbb{R}^2)$,
$$\|w\|_4^4\leq\frac{2}{\|Q\|_2^2}\|w\|_2^2\|\nabla w\|_2^2.$$
As a consequence, if  
$$\nL{w}\leq\nL{Q},$$ 
then the energy of $w$ is nonnegative.
 
Therefore on the one hand, 
$$0\leq E(e^{i\al\theta}v)$$ 
for every real number $\al$ and for all real function $\theta$, 
since $e^{i\alpha\theta}v$ is still a function of critical or
subcritical mass. 
On the other hand 
$$E(e^{i\al\theta}v)=
\frac{1}{2}\int\mo{ i\al\nabla\theta\, v+\nabla v}^2dx-
\frac{1}{4}\int\mo{v}^4dx$$ 
$$=\frac{\al^2}{2}\int\mo{v}^2\mo{\nabla\theta}^2dx-
\al\int\Im{(v\nabla\overline{v}\,)}\nabla\theta dx+E(v)$$ 
Thus the discriminant of the equation in $\al$ must be negative or null  
and we obtain the claimed Cauchy-Schwarz type inequality $(*)$.

\end{proof} 
 
\subsection{The concentration of the solution}
In this subsection we shall give a refined description of a critical
mass blowing-up solution $u$ of $(S_\Omega)$, by following the
approach of Weinstein in \cite{We33}.
 
In order to deal with real functions, we shall analyze the modulus of 
$u$. 
However, the same arguments below can be used to get the corresponding 
results on $u$ (see Remark \ref{Rim}).

One can write the convergence \eqref{limW}
$$u(t,x)=e^{-i\theta(t)}\la(t)(Q+R(t))(\la(t)(x-x(t))),$$
with $R$ a complex function such that 
$$\|R(t)\|_{\mathbb{H}^1(\R^2)}\underset{t\tend T}{\longrightarrow}0.$$
Since the modulus is a continuous function on $\mathbb{H}^1(\R^2)$ 
(\cite{MM3}), this 
 implies
$$|u(t,x)|=\la(t)|Q+R(t)|(\la(t)(x-x(t)))=
\la(t)(Q+\tilde{\tilde{R}}(t))(\la(t)(x-x(t))),$$
with $\tilde{\tilde{R}}$ a real function strongly converging to $0$ in 
$\mathbb{H}^1(\mathbb{R}^2)$. 

Let us set
$$\tilde{\la}(t)=
\frac{\nL{\nabla |u(t)|}}{\nL{\nabla Q}}.$$
By noticing that $|u(t)|$ is also of critical mass, its energy is 
nonnegative, and
\begin{equation}\label{difgrad}
0\leq \|\nabla u(t)\|_2^2-\|\nabla |u(t)|\|_2^2=2E(u)-2E(|u(t)|)\leq 
2E(u),
\end{equation}
which implies
$$(\la(t)+\tilde{\la}(t))(\la(t)-\tilde{\la}(t))=O(1).$$
Since $0\leq \tilde{\la}(t)\leq \la(t)$,
$$\la(t)-\tilde{\la}(t)=O\left(\frac{1}{\la(t)}\right),$$
and we have 
\begin{equation}\label{Rform}
|u(t,x)|=\tilde{\la}(t)(Q+\tilde{R}(t))(\tilde{\la}(t)(x-x(t))),
\end{equation}
with $\tilde{R}(t)$ a real function such that 
$$\|\tilde{R}(t)\|_{\mathbb{H}^1(\R^2)}
\underset{t\tend T}{\longrightarrow}0.$$

\begin{prop}\label{dec}
The remainder term $\tilde{R}$ has the decay 
\begin{equation}\label{Rdecay}
\|\tilde{R}(t)\|_{\mathbb{H}^1}\leq \frac{\tilde{C}}{\tilde{\la}(t)}
\leq \frac{C}{\la(t)}.
\end{equation}
\end{prop}

The proof follows Merle's one in \cite{Me13}. 
However, for the sake of completness, we give in the Appendix 
a proof by a slightly different method.

\bigskip

Finally, let us give the following property of decay of the solution.

\begin{lemma}\label{boundgrad}
Let $u$ be a critical mass solution of $(S_\Omega)$, 
blowing up at the finite time $T$, at one point $x_0\in\Omega$, which 
means that the concentration parameter $x(t)$ converges to $x_0$. 
Then, the gradient of $u(t)$ restricted outside any neighborhood $V$ of 
$x_0$ satisfies
$$\sup_{t\in[0,T[}\int_{^cV}|\nabla u(t)|^2dx<\infty.$$
\end{lemma}

\begin{proof}

The inequality \eqref{difgrad} implies
$$\sup_{t\in[0,T[}\int_{^cV}|\nabla u(t)|^2dx\leq
2E(u)+\sup_{t\in[0,T[}\int_{^cV}|\nabla |u(t)||^2dx.$$
By using \eqref{Rform},
$$\int_{^cV}|\nabla |u(t)||^2dx=
\tilde{\la}^2(t)\int_{^c\tilde{\la}(t)(V-x(t))}
|\nabla Q+\nabla\tilde{R}(t)|^2.$$
Since $x(t)$ converges to $x_0$ and $Q$ is exponentially decreasing, 
$$\tilde{\la}^2(t)\int_{^c\tilde{\la}(t)(V-x(t))}|\nabla Q|^2=o(1).$$
Then it follows that 
$$\int_{^cV}|\nabla |u(t)||^2dx\ls 
\tilde{\la}^2(t)\int_{^c\tilde{\la}(t)(V-x(t))}|\nabla \tilde{R}|^2,$$
and the decay \eqref{Rdecay} of $R$ implies 
$$\int_{^cV}|\nabla |u(t)||^2dx=O(1),$$
so the lemma is proved.
\end{proof}
\begin{remark}
Another proof of this lemma can be done by using the approach of Merle 
in \cite{Me3}. 
However, we shall need the full strength of Proposition \ref{dec} later 
in \S3.3 and \S3.4.
\end{remark}

\bigskip

\subsection{Derivatives of virial type functions}
Let $u$ be a solution of $(S_\Omega)$ and let $h$ be a 
$\mathcal{C}^\infty(\R^2)$ function with bounded first and second 
derivatives. 
Then, by using the fact that $u$ satisfies $(S_\Omega)$, we obtain
$$\partial_t\int_\Omega|u(t)|^2h dx=
2\int_\Omega \Re \left(u(t)\overline{u}_t(t)\right)h dx=
2\int_\Omega \Im \left(u(t)\Delta\overline{u}(t)\right)h dx.$$
Since $u$ cancels on the boundary of $\Omega$, by integration by parts 
\begin{equation}\label{g'}
\partial_t\int_\Omega|u(t)|^2h dx=
-2\int_\Omega \Im \left(u(t)\nabla\overline{u}(t)\right)\nabla h dx.
\end{equation}
By using again the equation $(S_\Omega)$
$$\DP{t}^2\int_\Omega\mo{u}^2h=
-2\int_\Omega \Re \left((\Delta u+|u|^2u)\nabla\overline{u}\right)
\nabla h+
2\int_\Omega \Re \left(u\nabla(\Delta \overline{u}+|u|^2\overline{u})
\right)
\nabla h$$
$$=\int_\Omega -
|u|^2\Delta^2h-
|u|^4\Delta h+
2|\nabla u|^2\Delta h-
4\Re\left(\Delta u\nabla\overline{u}\right)\nabla h.$$
It follows that
$$\DP{t}^2\int_\Omega\mo{u}^2h=
\int_\Omega -
|u|^2\Delta^2h-
|u|^4\Delta h+
4\Re \sum_{i,j}\partial_iu\partial_j\overline{u}\partial_{ij}h-
2\int_{\partial\Omega}\left|\frac{\partial u}{\partial \nu}\right|^2
\frac{\partial h}{\partial \nu}d\sigma.$$
Therefore, by making the energy of the solution appear, we have the 
following identity.
\begin{lemma}\label{g''=}
For a solution $u$ of $(S_\Omega)$ and a $\mathcal{C}^\infty(\R^2)$ 
function $h$ with bounded derivatives $\partial_{ij}h$ and $\Delta^2h$, 
we have 
$$\DP{t}^2\int_\Omega\mo{u}^2h=
16E(u)-\int_\Omega(2\mo{\nabla u}^2-\mo{u}^4)(4-\Delta h)-
\int_\Omega 
|u|^2\Delta^2h$$
$$+\int_{\Omega}\left(4\Re \sum_{i,j}\partial_iu\partial_j\overline{u}\partial_{ij}h-
2|\nabla u|^2\Delta h\right)-
2\int_{\partial\Omega}\left|\frac{\partial u}{\partial \nu}\right|^2
\frac{\partial h}{\partial \nu}d\sigma.$$
\end{lemma}
\begin{corollary}\label{g''}
For a solution $u$ of $(S_\Omega)$ and a $\mathcal{C}^\infty(\R^2)$ 
function $h$ equal to $|x|^2$ on $B(0,R)$, with bounded derivatives 
$\partial_{ij}h$ and $\Delta^2 h$, we have the estimate
$$
\left|\partial_t^2\int_\Omega |u(t)|^2hdx-16E(u)\right|\leq
C\int_{\{|x|\geq R\}\cap\Omega}(|u(t)|^2+|\nabla u(t)|^2)dx$$
$$+\int_{\partial\Omega}
\left|\frac{\partial u(t)}{\partial \nu}\right|^2\left|
\frac{\partial h}{\partial \nu}\right| d\sigma.$$
\end{corollary}

\section{The blow-up rate on bounded plane domains} 
\subsection{The convergence of the concentration points $x(t)$} 
\begin{lemma}
Let $\Omega$ be a bounded domain and let $u$ be a critical mass
solution of $(S_\Omega)$, blowing up at the finite time $T$. 
Then the concentration parameter $x(t)$ has a limit at the time $T$. 
\end{lemma}

\begin{proof}

From \eqref{limx} it follows that for a test function $\psi$, 
$$\int_{\Omega-x(t)}\mo{u(t,x+x(t))}^2\psi(x)dx
\underset{t\tend T}{\longrightarrow}\nL{Q}^2\psi(0).$$ 
If $\psi$ is chosen such that $\psi(0)\neq0$ then, since  
the set $\Omega$ is bounded, it follows that 
\begin{equation}\label{unifboundx}
\underset{t\tend T}{\limsup}\mo{x(t)}<\infty.
\end{equation}

The first order momentum  
$$f(t)=\int_\Omega\mo{u(t,x)}^2xdx,$$ 
stays finite in time since $\Omega$ is bounded and $u$ 
conserves its mass. 
By using the formula \eqref{g'} for vector-valued functions $h$, one can 
calculate the derivative 
$$f'(t)=-2\int_\Omega\Im{(u(t)\nabla\overline{u}(t)\,)}dx.$$ 
The inequality $(*)$ in the special case $\theta_i(x)=x_i$ implies 
that this derivative is bounded in time 
$$\mo{f'(t)}^2\leq
4\sum_{i\in\{1,2\}}
\left\vert\int_\Omega\Im{(u(t)\nabla\overline{u}(t)\,)}\nabla\theta_idx
\right\vert^2\leq 
16E(u)\|u\|_2^2.$$ 
Therefore $f$ admits a limit at the time $T$. 
Let us define $x_0$ by
$$f(T)=x_0\nL{Q}^2.$$ 
Using the convergence \eqref{limx} and \eqref{unifboundx} which implies
that $\Omega-x(t)$ is a uniformly bounded set, one has 
$$f(t)-x(t)\nL{Q}^2=
\int_{\Omega-x(t)}\mo{u(t,x+x(t))}^2xdx 
\underset{t\tend T}{\longrightarrow}0.$$

Therefore the point $x_0$ is the limit of $x(t)$, and the square of the 
solution behaves like a Dirac function 
\begin{equation}\label{lim}
\mo{u(t,\cdot)}^2
\underset{t\tend T}{\longrightarrow}\nL{Q}^2\,\delta_{x_0}.
\end{equation} 

\end{proof}
In the following, we shall suppose, up to a translation, 
that the solution blows up at the point $0\in\overline{\Omega}$.

\subsection{Lower bound for the blow-up rate} \label{gprim0}
The derivative in time of the the virial of the solution $u$, 
$$g(t)=\int_\Omega\mo{u(t,x)}^2\mo{x}^2dx,$$
can be calculated with the formula \eqref{g'} with $h(x)=|x|^2$, and
$$g'(t)=-4\int_\Omega\Im{(u(t)\nabla\overline{u}(t)\,)}xdx.$$ 
Therefore the inequality $(*)$ in the case $\theta(x)=\mo{x}^2$ implies 
that 
$$|g'(t)|\leq 4\sqrt{2E(u)g(t)}.$$ 
The concentration result \eqref{lim} of the former subsection gives 
$$g(T)=0,$$ 
and one can now write 
$$\sqrt{g(t)}=
-\int_t^T \frac{g'(\tau)}{2\sqrt{g(\tau)}}d\tau\leq
\int_t^T 2\sqrt{2E(u)}=
2\sqrt{2E(u)}(T-t),$$ 
and obtain 
$$g(t)\leq 8E(u)(T-t)^2.$$ 
Then the uncertainty principle 
$$\left(\int_{\R^2}\mo{u}^2\right)^2\leq
\left(\int_{\R^2}\mo{u}^2\mo{x}^2\right)
\left(\int_{\R^2}\mo{\nabla u}^2\right)$$ 
gives us a lower bound of the blow-up rate
$$\frac{\nL{Q}^2}{2\sqrt{2E(u)}(T-t)}\leq
\nL{\nabla u(t)},$$ 
so the first assertion of Theorem \ref{art31} is proved.
\subsection{Equivalence between the virial and the blow-up rate}

By using \eqref{Rform},
$$\int_\Omega |u(t)|^2|x-x(t)|^2dx=
\frac{1}{\tilde{\la}^2(t)}\int_{\tilde{\la}(t)(\Omega-x(t))}
(Q+\tilde{R}(t))^2|x|^2dx.$$
Since $x(t)$ tends to $0$ and $Q$ is exponentially decreasing,
$$\frac{1}{\tilde{\la}^2}\int_{\tilde{\la}(t)(\Omega-x(t))}Q^2|x|^2dx=
O\left(\frac{1}{\tilde{\la}^2(t)}\right),$$
so 
$$\int_\Omega |u(t)|^2|x-x(t)|^2dx\ls 
\frac{1}{\tilde{\la}^2(t)}\int_{\tilde{\la}(t)(\Omega-x(t))}
\tilde{R}^2(t)|x|^2dx+\frac{1}{\tilde{\la}^2(t)}.$$ 
The domain $\Omega$ is considered bounded, so one can write
$$\int_\Omega |u(t)|^2|x-x(t)|^2dx\ls 
\int\tilde{R}^2(t)dx+\frac{1}{\tilde{\la}^2(t)},$$ 
and by using the decay \eqref{Rdecay} of $\tilde{R}$, we obtain
$$\int_\Omega |u(t)|^2|x-x(t)|^2dx\ls \frac{1}{\la^2(t)}.$$
As we did in the previous subsection, by the uncertainy principle for 
$u(t,x+x(t))$,
$$\|u\|_2^4\ls \la^2(t)\int |u(t)|^2|x-x(t)|^2dx,$$
and so the first assertion of Proposition \ref{equivir} follows,
$$\int_\Omega |u(t)|^2|x-x(t)|^2dx\approx \frac{1}{\la^2(t)}.$$

\bigskip

\subsection{A differentiable choice for $x(t)$}
Let us set 
$$y(t)=\frac{\int |u(t)|^2xdx}{\|Q\|^2}.$$
By using the conservation of the mass, which is critical,
$$x(t)-y(t)=\frac{1}{\|Q\|_2^2}\int |u(t)|^2(x-x(t))dx.$$
Then by \eqref{Rform} one has
$$x(t)-y(t)=
\frac{1}{\tilde{\la}(t)\|Q\|_2^2}
\int_{\tilde{\la}(t)(\Omega-x(t))}(Q+R(t))^2 xdx.$$
Therefore, by the same arguments as in the previous subsection, and 
the by using the fact that since $Q$ is radially symmetric,
$$\int Q^2(x)xdx=0,$$
then
$$|x(t)-y(t)|\leq \frac{C}{\tilde{\la}^2(t)}.$$ 
If we define $S$ by
$$|u(t,x)|=\tilde{\la}(Q+\tilde{R}(t))(\tilde{\la}(x-x(t)))=
\tilde{\la}(Q+S(t))(\tilde{\la}(x-y(t))),$$
one has
$$\|S(t)\|_{\mathbb{H}^1}\leq 2\|\tilde{R}(t)\|_{\mathbb{H}^1}+
\|Q(\cdot+\la(x(t)-y(t))-Q(\cdot)\|_{\mathbb{H}^1}.$$
The decay of the difference between $x(t)$ and $y(t)$, together 
with \eqref{Rdecay}, implies
$$\|S(t)\|_{\mathbb{H}^1}\leq \frac{C}{\tilde{\la}(t)}.$$ 
So, by changing $x(t)$ into 
$$\frac{\int |u(t)|^2xdx}{\|Q\|^2},$$
we have the convergence corresponding to \eqref{Rform}
$$|u(t,x)|=\tilde{\la}(t)(Q+S(t))(\tilde{\la}(t)(x-y(t))),$$
with $S$ decreasing in $\mathbb{H}^1$ as does $R$, and so the second 
assertion of Proposition \ref{equivir} follows.

The interest of this choice of the concentration parameter is that 
$y(t)$ is a differentiable function, and, moreover, in the radial case 
we obtain Corollary \ref{radial}.
\section{The blow-up rate on unbounded plane domains} \label{gloc}

Consider now the equation $(S_\Omega)$ on an unbounded domain of the plane or 
on a surface. Let $u$ be a critical mass solution that blows up in
an interior point $x_0$ of $\Omega$, that is
$$x(t)\underset{t\tend T}{\longrightarrow}x_0.$$ 
Modulo a translation, we can suppose that $x_0$ is zero 
and so, 
$$\mo{u(t,x)}^2
\underset{t\tend T}{\longrightarrow}\nL{Q}^2\,\delta_0.$$

Let $\phi$ be a $\mathcal{C}^\infty_0$ function, equal to $1$ on 
$B(0,R)$. 
Let us introduce the localized virial of the solution 
$$g_\phi(t)=\int\mo{u(t,x)}^2\phi^2(x)\mo{x}^2dx.$$ 
Then, using \eqref{g'} with $h(x)=\phi^2(x)\mo{x}^2$, one has
$$g_\phi'(t)=
-2\int\Im(u(t)\nabla\overline{u}(t))\nabla(\phi^2\mo{x}^2)dx.$$ 
The inequality $(*)$ with $\theta(x)=\phi^2(x)\mo{x}^2$ gives us
$$\mo{g_\phi'(t)}^2\leq
8E(u)\int\mo{u}^2\mo{\nabla(\phi^2\mo{x}^2)}^2dx$$
Since $\nabla(\phi^2\mo{x}^2)$ is a $\mathcal{C}^\infty_o(\R^2)$ 
function cancelling at $0$, and since the square of $\mo{u}$ behaves 
like a Dirac distribution, it follows that
$$g_\phi'(T)=0.$$ 
Then, as in the former section, and using the existence of a positive 
constant $C$ such that 
$$|\nabla (\phi^2|x|^2)|^2\leq C\psi^2|x|^2,$$
one has 
$$g_\phi(t)\ls (T-t)^2.$$ 
The uncertainty principle reads 
$$\left(\int\mo{u}^2\phi^2 dx\right)^2\leq 
\left(\int\mo{u}^2\phi^2\mo{x}^2dx\right)
\left(\int\mo{\nabla (u\phi)}^2dx\right).$$
By integrating by parts the last term and by using the fact that $\phi$ 
is equal to $1$ on $B(0,R)$, it follows that
$$\left(\int_{B(0,R)}\mo{u(t)}^2\right)^2\leq 
g_\phi(t)
\left(\int\mo{\nabla u}^2\phi^2dx-\int|u|^2\phi\Delta\phi dx\right).$$
Since $\phi$ is a $\mathcal{C}^\infty_0$ function,
$$\left(\int_{B(0,R)}\mo{u(t)}^2\right)^2\leq 
g_\phi(t)
\left(C\int\mo{\nabla u}^2dx-\int|u|^2\phi\Delta\phi dx\right).$$
On the one hand the $\mathbb{L}^2$ norm of $u$ is conserved. 
On the other hand, the behavior of $|u|^2$ as a Dirac distribution 
implies that the norm of 
its restriction outside a neighborhood of zero tends to $0$ in time.
So we have
$$\left\{\begin{array}{c}
\int_{B(0,R)}\mo{u(t)}^2=O(1),\\
\\
\int|u(t)|^2\phi\Delta\phi dx=o(1),
\end{array}\right.$$
and since $g_\phi$ is bounded in time,
$$1\ls\sqrt{g_\phi (t)}\nL{\nabla u(t)}$$
Then the decay of $g_\phi$ gives us the lower bound of the blow-up 
speed 
$$\frac{1}{T-t}\ls\nL{\nabla u(t)}.$$


\section{Blow-up on the boundary} 

\subsection{Necessary condition for blow-up on the boundary} 
Let us first introduce a notion of limit of sets, as in \cite{IPG3}.
\begin{definition}\label{def}
A sequence of open sets $M_m$ is said to tend to an open set $M$ of 
$\R^2$ if the following conditions are satisfied.\\
i) For all compact $K\subset M$, there exists $n_K\in\mathbb{N}$, such that 
for all $n\geq n_K$, $K\subset M_n$.\\
ii) For all compact $K\subset ^c\overline{M}$, there exists 
$n_K\in\mathbb{N}$, such that for all $n\geq n_K$, 
$K\subset ^c\overline{M_n}$.
\end{definition}

Let us suppose that there exists an explosive solution $u$ of 
the equation $(S_\Omega)$ at $0\in\partial\Omega$. 
The convergence \eqref{limW} implies that 
$$\la(t)(\Omega-x(t))\underset{t\tend T}{\longrightarrow}\mathbb{R}^2.$$ 
As in \cite{IPG3}, the limit set depends on the position of $x(t)$ with 
respect to the boundary of $\Omega$. 
If there is a positive number $C$ such that for all $t$ 
$$\la(t)d(x(t),\partial\Omega)\leq C,$$ 
then $\la(t)(\Omega-x(t))$ tends to a half-plane 
and blow-up cannot occur. Also, if 
$$\la(t)d(x(t),\partial\Omega)
\underset{t\tend T}{\longrightarrow}\infty,$$ 
and $x(t)$ is not in $\Omega$, then, by Definition \ref{def}, 
$\la(t)(\Omega-x(t))$ tends to the empty set. 
Therefore the only possibility to have explosion on the boundary is 
that $x(t)\in\Omega$ and 
$$\la(t)d(x(t),\partial\Omega)
\underset{t\tend T}{\longrightarrow}\infty.$$ 
In particular, since $0$ is on the boundary, 
\begin{equation}\label{bound}
\mo{\la(t)x(t)}\underset{t\tend T}{\longrightarrow}\infty.
\end{equation} 
 
We have 
$$\mo{x(t)}^2\int_{B(x(t),\frac{C}{\la(t)})}\mo{u}^2\leq
2\int_{B(x(t),\frac{C}{\la(t)})}\mo{u}^2\mo{x-x(t)}^2+
2\int_{B(x(t),\frac{C}{\la(t)})}\mo{u}^2\mo{x}^2.$$
On the one hand, by using the Weinstein relation \eqref{limW}, one has 
$$\mo{x(t)}^2\int_{B(x(t),\frac{C}{\la(t)})}\mo{u}^2\approx 
\mo{x(t)}^2.$$ 
On the other hand, using again \eqref{limW}, 
$$\int_{B(x(t),\frac{C}{\la(t)})}\mo{u}^2\mo{x-x(t)}^2\ls
\frac{1}{\la(t)^2}.$$ 
In view of \eqref{bound}, these two facts imply 
$$\mo{x(t)}^2\ls\int_{B(x(t),\frac{C}{\la(t)})}\mo{u}^2\mo{x}^2\ls 
g_\psi,$$
where $g_\psi$ is the localized virial function defined in \S\ref{gloc}. 
In the same section it was proved that 
$$g_\psi\ls (T-t)^2,$$ 
so it follows that
$$\mo{x(t)}\ls T-t.$$
By using again \eqref{bound}, 
$$\frac{1}{T-t}\ll\la(t),$$ and the second assertion of 
Theorem \ref{art31} is proved.

\subsection{Results of non-explosion}

From now on we assume that 
$\Omega$ be a half plane whose boundary contains $0$ or a 
plane sector with corner $0$.  
Suppose there exists an explosive solution $u$ of critical mass such 
that $u$ behaves like a Dirac mass at $0$.

For a radial function $f\in\mathcal{C}^\infty(\R^2)$, 
the result of Lemma \ref{g''=} becomes
$$\DP{t}^2\int_\Omega\mo{u}^2f=
16E(u)-\int_\Omega(2\mo{\nabla u}^2-\mo{u}^4)(4-\Delta f)-
\int_\Omega 
|u|^2\Delta^2f$$
\begin{equation}\label{g''b=}
+\int_{\Omega}4\Re \sum_{i,j}\partial_iu\partial_j\overline{u}\partial_{ij}f-
2|\nabla u|^2\Delta f,
\end{equation}
since from the choice of $\Omega$ 
$$x.\nu=0 \mbox{ on }\partial\Omega.$$ 
It follows that for a radial function $f\in\mathcal{C}^\infty(\R^2)$, 
equal to $|x|^2$ on $B(0,R)$, with bounded derivatives $\partial_{i,j}f$ and 
$\Delta^2 f$, the estimate of Corollary \ref{g''} becomes
\begin{equation}\label{g''b}
\left|\partial_t^2\int_\Omega |u(t)|^2f-16E(u)\right|\leq
C\int_{\{|x|\geq R\}\cap\Omega}|u(t)|^2+|\nabla u(t)|^2.
\end{equation}

Arguing as in \cite{Me3}, we obtain the following lemmas.

\begin{lemma}\label{varfin}
The initial condition is of finite variance
$$\int_\Omega\mo{u_0}^2\mo{x}^2dx<\infty.$$ 
\end{lemma}

\begin{proof}
Let us consider $\psi$ a $\mathcal{C}_0^\infty(\R)$ positive radial 
function which is equal to $|x|^2$ on $B(0,1)$. 
Notice that
$$\mo{\nabla\psi}^2\leq C\psi.$$
For all entire $n$, we introduce the localized virial functions
$$g_n(t)=\int_\Omega\mo{u(t)}^2\psi_ndx,$$
where
$$\psi_n(x)=n^2\psi\left(\frac{x}{n}\right).$$
The Taylor formula in zero for the function $g_n(t)$ gives us
$$\mo{g_n(t)-g_n(0)}\leq t\mo{g_n'(0)}+C\sup_{t}\mo{g_n''(t)}.$$
Since $\psi_n$ are equal to $|x|^2$ on $B(0,1)$, and the derivatives 
$\partial_{ij}\psi_n$ and $\Delta^2\psi_n$ are uniformly bounded, we 
can estimate by \eqref{g''b} 
$$\left\vert g_n''(t)-16E(u)\right\vert\leq 
C\int_{\mo{x}>1}(\mo{u(t)}^2+\mo{\nabla u(t)}^2)dx.$$ 
Then, in view of Lemma \ref{boundgrad}, the quantity $g_n''(t)$ is bounded 
uniformly on $n$. 
So we have
$$\mo{g_n(t)-g_n(0)}\leq T\mo{g_n'(0)}+C.$$
By using the inequality $(*)$, 
$$\mo{g_n'(t)}=
\left\vert\int_\Omega\Im{(u(t)\nabla\overline{u}(t)\,)}\nabla\psi_n
\right\vert\leq 
\left(2E(u)\int_\Omega\mo{u(t)}^2\mo{\nabla\psi_n}^2
\right)^{\frac{1}{2}}.$$ 
The choice of $\psi_n$ gives us
$$|\nabla\psi_n|^2\leq \psi_n,$$
and it follows that
$$\mo{g_n'(t)}\leq 
C\left(\int_\Omega\mo{u(t)}^2\psi_n\right)^{\frac{1}{2}}=
C\sqrt{g_n(t)}.$$ 
Therefore 
$$g_n(0)-2C\sqrt{g_n(0)}-C\leq g_n(t)$$ 
for any time $t$, and 
$$g_n(0)-2C\sqrt{g_n(0)}-C\leq \underset{t\tend T}{\lim}g_n(t).$$ 
The concentration of the solution as a Dirac distribution 
implies that for fixed $n$ 
$$\underset{t\tend T}{\lim}g_n(t)=0,$$ 
and therefore
$$\underset{n\tend \infty}{\lim}(g_n(0)-2C\sqrt{g_n(0)}-C)\leq 0.$$
As a consequence, $g_n(0)$ is bounded as $n$ tends to infinity. 
Since the supports of $\psi_n$ cover $\Omega$ when $n$ tends to infinity, 
it follows that the initial condition is of finite variance
$$\int_\Omega\mo{u_0}^2\mo{x}^2dx<\infty.$$ 
 
\end{proof}

\begin{remark}
When $\Omega$ is a bounded domain, and $u$ is a critical mass function 
blowing up at a point of $\Omega$ or of its boundary, 
it is easy to see that the initial condition is of finite variance
$$\int_\Omega|u|^2|x|^2dx\leq C\|u\|_2^2<\infty.$$
\end{remark}


\begin{lemma}\label{limg} 
The limit in time of the virial function is
$$g(T)=0.$$
\end{lemma}

\begin{proof}
Let us consider a $\mathcal{C}^\infty$ positive function $\phi$ which 
is null on $B(0,1)$ and verifies
$$\frac{\mo{x}}{2}\leq\phi(x)\leq \mo{x},$$
on $^cB(0,2)$. 
Suppose also that the derivatives 
$\partial_{ij}\psi$ and $\Delta^2\psi$ are bounded.
We denote 
$$\phi_n(x)=n\,\phi\left(\frac{x}{n}\right),$$ 
so $\phi_n$ are supported on $^cB(0,n)$ and verify
$$\frac{|x|}{2}\leq\phi_n(x)\leq |x|$$
on $^cB(0,2n)$.

Taylor's formula together with \eqref{g'} and the estimate 
\eqref{g''b} gives us 
$$\int\mo{u(t)}^2\phi_n\leq \int_{\mo{x}>n}\mo{u_0}^2\mo{x}^2+
T\left\vert\int_{\mo{x}>n}
\Im{(u_0\nabla\overline{u_0}\,)}\nabla\phi_n^2\right\vert$$
\begin{equation}\label{g(T)}
+C(T-t)^2+C\int_0^T(T-\tau)\int_{\mo{x}>n}
(\mo{\nabla u(\tau)}^2+\mo{u(\tau)}^2)d\tau.
\end{equation}
The Lemma \ref{varfin} ensures us that the initial data is of finite 
variance, therefore
$$\int_{\mo{x}>n}\mo{u_0}^2\mo{x}^2
\underset{n\tend\infty}{\longrightarrow}0,$$
Also, by Lemma \ref{boundgrad},
$$\left\vert
\int_{\mo{x}>n}\Im{(u_0\nabla\overline{u_0}\,)}\nabla\phi_n^2\right\vert
\underset{n\tend\infty}{\longrightarrow}0.$$
Then, using again Lemma \ref{boundgrad} and the conservation of the 
mass, for all $\tau$ and for all $n$ there exist a positive constant 
$C$ such that 
$$\int_{\mo{x}>n}(\mo{\nabla u(\tau)}^2+\mo{u(\tau)}^2)\leq C.$$
One also has, for every $\tau$,
$$\int_{\mo{x}>n}(\mo{\nabla u(\tau)}^2+\mo{u(\tau)}^2)
\underset{n\tend\infty}{\longrightarrow}0.$$
Then by the dominated convergence theorem
$$\int_0^T(T-\tau)\int_{\mo{x}>n}
(\mo{\nabla u(\tau)}^2+\mo{u(\tau)}^2)d\tau
\underset{n\tend\infty}{\longrightarrow}0.$$
Therefore it follows from \eqref{g(T)} that for all $t$,
$$\int\mo{u(t)}^2\phi_n^2\leq\epsilon (n)+C(T-t)^2,$$
with
$$\epsilon(n)\underset{n\tend\infty}{\longrightarrow}0.$$
On the one hand, in view of the choice of $\phi_n$, this gives us
$$\int_{\mo{x}>2n}\mo{u(t)}^2\mo{x}^2
\leq 2\epsilon(n)+C(T-t)^2.$$
On the other hand, for fixed $n$, the concentration of the 
solution as a Dirac distribution implies
$$\underset{t\tend T}{\lim}\int_{\mo{x}<2n}\mo{u(t)}^2\mo{x}^2=0.$$
Therefore, for all $n$
$$\underset{t\tend T}{\lim}\int\mo{u(t)}^2\mo{x}^2
\leq\epsilon(n).$$
By letting $n$ to tend to infinity one has
$$\underset{t\tend T}{\lim}\int\mo{u(t)}^2\mo{x}^2=0,$$
that is
$$g(T)=0,$$
and the Lemma \ref{limg} is proved.
\end{proof} 

This lemma and the same arguments as in \S\ref{gprim0} give us also
$$g'(T)=0.$$ 
By using the formula \eqref{g''b=} with $f(x)=|x|^2$, 
the second derivative of the virial is exactly 
$$g''(t)=16E(u).$$
Then it follows that 
$$g(t)=8E(u)(T-t)^2,$$ 
and by the same calculation as in \S\ref{CS} 
$$E(e^{i\frac{\mo{x}^2}{4(T-t)}}u(t,x))=
E(u)+\frac{1}{4(T-t)}g'(t)+\frac{1}{16(T-t)^2}g(t)=0.$$ 
For fixed $t$, by the variational characterization of the ground state 
$Q$, there exists real numbers $\theta$ and $\omega$ such that 
$$u(t,x)=e^{-i\frac{\mo{x}^2}{4(T-t)}}e^{i\theta}\omega \,Q(\omega 
(x-x_0))$$
for some $x_0\in\mathbb{R}^2$ (\cite{Ca03}). 
This means that the support of $u$ is the entire $\mathbb{R}^2$ that is 
a contradiction, and the proof of Theorem \ref{art32} is complete.
 

\section{Appendix}
In this Appendix we give a proof for Proposition \ref{dec}.
Let us recall the notations of \S2.2. 
We have defined
$$\tilde{\la}(t)=
\frac{\nL{\nabla |u(t)|}}{\nL{\nabla Q}},$$
and the solution $u$ was written \eqref{Rform}
$$|u(t,x)|=\tilde{\la}(t)(Q+\tilde{R}(t))(\tilde{\la}(t)(x-x(t))),$$
with $\tilde{R}(t)$ a real function such that 
$$\|\tilde{R}(t)\|_{\mathbb{H}^1(\R^2)}
\underset{t\tend T}{\longrightarrow}0.$$
We shall prove in the following the decay \eqref{Rdecay} 
asserted in Proposition \ref{dec} 
$$\|\tilde{R}(t)\|_{\mathbb{H}^1}\leq \frac{\tilde{C}}{\tilde{\la}(t)}
\leq \frac{C}{\la(t)}.$$

\bigskip

The fact that $u$ is of critical mass gives us
\begin{equation}\label{massR}
\int \tilde{R}^2=-2\int Q\tilde{R},
\end{equation}
and the choice of $\tilde{\la}$ implies
$$\int |\nabla \tilde{R}|^2=-2\int\nabla Q\nabla \tilde{R}.$$
Let us calculate the energy of $|u|$,
$$\frac{2E(|u|)}{\tilde{\la}^2}=\int|\nabla Q+\nabla \tilde{R}|^2-
\frac{1}{2}\int
(Q+\tilde{R})^4.$$
The energy of $Q$ is zero, so
$$\frac{2E(|u|)}{\tilde{\la}^2}=\int|\nabla \tilde{R}|^2+2\nabla Q\nabla 
\tilde{R}-\frac{\tilde{R}^4}{2}-
2Q\tilde{R}^3-3Q^2\tilde{R}^2-2Q^3\tilde{R}.$$
The ground state $Q$ verifies the equation
$$\Delta Q+Q^3=Q,$$
and therefore, by using the relation \eqref{massR} on $\tilde{R}$,
$$\int2\nabla Q\nabla \tilde{R}-2Q^3\tilde{R}=-2\int Q\tilde{R}=
\int \tilde{R}^2.$$
So finally
$$<L \tilde{R},\tilde{R}>=\frac{2E(|u|)}{\tilde{\la}^2}+
\frac{1}{2}\int \tilde{R}^4+\int 2Q\tilde{R}^3,$$
where $L$ is the operator 
$$L=-\Delta+(1-3Q^2).$$
Since $\tilde{R}$ tends to $0$ in $\mathbb{H}^1$, by using the Sobolev
embeddings, the cubic and quadratic terms in $\tilde{R}$ are negligible with
respect to the $\mathbb{H}^1$ norm of $\tilde{R}$.
Also, the energy of $|u|$ is bounded by the constant energy of $u$,
so for having \eqref{Rdecay} it is 
sufficient to prove the existence of a positive constant $\delta$ such
that for $t$ close enough to $T$
$$\delta\|\tilde{R}(t)\|_{\mathbb{H}^1}^2\leq<L \tilde{R}(t),\tilde{R}(t)>.$$

\begin{remark}\label{Rim}
The initial complex function $R$ can be analyzed in the same manner, 
and one has 
$$<L_{-}\Im R,\Im R>+<L\Re R,\Re R>\leq 
\frac{2E(u)}{\la^2}+\frac{1}{2}\int |R|^4+\int 2Q|R|^3,$$ 
where $L_{-}$ is the operator
$$L_{-}=-\Delta +(1-Q).$$
This operator is non-negative and its kernel is spanned by $Q$. 
So once the decay \eqref{Rdecay} is obtained, by decomposing $\Im R$ with 
respect to $Q$, we also have 
$$\|R(t)\|_{\mathbb{H}^1}\leq \frac{C}{\la(t)}.$$
\end{remark}

Following the ideas of Weinstein in \cite{We33}, we shall look for the 
nature of the negative eigenvalues of $L$.
\begin{lemma}\label{2vp}
The second eigenvalue of $L$ is $0$.
\end{lemma}
\begin{proof}
Let us consider the functional
$$J(f)=\frac{\|u\|_2^2\|\nabla u\|_2^2}{\|u\|_4^4},$$
which is minimized by $Q$ (see the introduction). 
Then, for a test function $f$,
$$\partial_{\epsilon}^2 J(Q+\epsilon f)|_{\epsilon =0}\geq 0.$$
By explicitly calculating this second derivative and using \eqref{normQ} in 
the calculus, one has 
$$2\|Q\|_2^2<Lf,f>\geq -8<Q,f><\nabla Q,\nabla f>.$$
If we take $f$ to be orthogonal to $Q$, then
$$<Lf,f>\geq 0,$$
and by the Min-Max Principle (\cite{RS3}), 
the second eigenvalue of $L$ is non-negative. 
By noticing that the two partial derivatives of $Q$ verify 
$$L\partial_iQ=0,$$
we obtain that $0$ is an eigenvalue of $L$ of order grater than one, so 
the first eigenvalue is negative. 
Therefore the second eigenvalue of $L$ is $0$.
\end{proof}

We shall use the following theorem.
\begin{theorem}\label{mihai} (Maris \cite{M3}). 
Let $g\in\mathcal{C}^1([0,\infty))$, with $g(0)=0$, $g'(0)>0$ and 
$|g'(s)-g'(0)|\leq C|s|^\alpha$, for small $s$ and some $C,\alpha >0$. 
Let $a_0=\sup\{a>0|g(s)>0,\forall s\in(0,a)\}$, and let 
$u_0$ be a ground state of the operator
$$-\Delta u+g(u).$$
We define
$$I(u,\la)=\la ug'(u)-(\la+2)g(u),$$
and we will make the following assumptions : $a_0<u_0(0)$ and 
there exists a continuous 
function $\la : (a_0,u_0(0)]\tend(0,\infty)$ such that for any $U\in 
(a_0,u_0(0)]$ we have 
$$\left\{\begin{array}{c}
I(u,\la(U))\leq 0, \forall u\in[0,U],\\ 
I(u,\la(U))\geq 0, \forall u\in[U,u_0(0)].
\end{array}\right.$$
Then 
$$Ker (-\Delta+g'(u_0))=\{\partial_1u_0,\partial_2u_0\}.$$
\end{theorem}
Next we show that the operator $L$ satisfies the hypothesis of the 
theorem.
\begin{lemma}\label{ker}
The kernel of $L$ has dimension $2$.
\end{lemma}
\begin{proof}
In we take function $g$ to be
$$g(s)=s-s^3,$$
then $a_0=1$, the ground state $u_0$ is $Q$,
$$-\Delta+g'(u_0)=L,$$
and
$$I(u,\la)=2u((1-\la)u^2-1).$$
Let us consider the integral of $g$,
$$G(s)=\frac{s^2}{2}-\frac{s^4}{4}.$$
By using the relation \eqref{normQ} between the 
$\mathbb{L}^2$ and the $\mathbb{L}^4$ norms of $Q$ 
$$\int G(Q(x))dx=0.$$
The positivity of $G(s)$ on $[0,\sqrt{2}[$ implies the existence  
of points $x$ such that $Q(x)>\sqrt{2}$, and in particular 
$Q(x)>1$. 
Let us recall that $Q$ is a radial positive decreasing 
function. It follows that $Q(0)>1$, and the first assumption of 
the theorem \ref{mihai} is satisfied. 
The second assumption is satisfied for the function 
$$\la(U)=1-\frac{1}{U^2},$$
and we can conclude that 
$$Ker L=\{\partial_1Q,\partial_2Q\}.$$
\end{proof}

We return now to the study of $\tilde{R}$. 
We impose a choice of $x(t)$ which will yield an orthogonality 
property of $\tilde{R}$. 
Since
$$\frac{1}{\tilde{\la}(t)}
|u|\left(t,\frac{x}{\tilde{\la}(t)}+x(t)\right)
\underset{t\tend T}{\longrightarrow}Q(x),$$
we can choose $x(t)$ such that the functional
$$I(z)=\left\|\frac{1}{\tilde{\la}(t)}
|u|\left(t,\frac{\cdot+z}{\tilde{\la}(t)}+x(t)\right)-Q(\cdot)
\right\|_{\mathbb{H}^1}^2$$
reaches its minimum for $z=0$. 
By using \eqref{Rform}, this implies that the derivative in $z$ of 
$$\partial_z \|(Q+\tilde{R}(t))(\cdot)-Q(\cdot-z)\|_{\mathbb{H}^1}^2,$$
must be zero at $z=0$. 
It follows that 
$$\int\tilde{R}\,\partial_iQ+
\int\nabla \tilde{R}\,\partial_i\nabla Q=0.$$
One can then integrate by parts and obtain
$$\int\tilde{R}\,\partial_iQ-\int\tilde{R}\,\partial_i\Delta Q=0.$$
By recalling that the ground state $Q$ verifies 
$$\Delta Q+Q^3=Q,$$
it follows that $\tilde{R}$ has the orthogonality property
\begin{equation}\label{ortho3}
<\partial_iQ^3,\tilde{R}>=0.
\end{equation}

Let us recall that for having the decay property \eqref{Rdecay} of 
$\tilde{R}$, it is sufficient to prove that the operator $L$ controls 
its $\mathbb{H}^1$ norm.

\begin{lemma}
There exist a positive constant $\delta$ such that for $t$ close
enough to $T$,
$$\delta\|\tilde{R}(t)\|_{\mathbb{H}^1}^2\leq<L \tilde{R}(t),\tilde{R}(t)>.$$
\end{lemma}

\begin{proof}
We denote by $R_\shortparallel$ the projection of $\tilde{R}$ on the 
space spanned by $Q$, and by $R_\bot$ the remainder term, orthogonal to 
$Q$. 
Since the operator $L$ is self-adjoint,
$$<L\tilde{R},\tilde{R}>=<LR_\shortparallel,R_\shortparallel>+
2<LR_\shortparallel,R_\bot>+<LR_\bot,R_\bot>.$$
The first term reads
$$<LR_\shortparallel,R_\shortparallel>=<LQ,Q>\frac{<Q,\tilde{R}>^2}{\|Q\|_2^4},$$
and by using \eqref{massR}
$$<LR_\shortparallel,R_\shortparallel>=C\|\tilde{R}\|_2^4.$$

The second term is
$$<LR_\shortparallel,R_\bot>=\frac{<Q,\tilde{R}>}{\|Q\|_2^2}<LQ,R_\bot>,$$
and since $LQ=-2Q^3$, by using the Cauchy-Schwarz inequality, 
$$<LR_\bot,R_\shortparallel>=-2\frac{<Q,\tilde{R}>}{\|Q\|_2^2}<Q^3,R_\bot>\leq 
C\|\tilde{R}\|_2^3.$$

Now we have to estimate the third term. 
Let us notice that the orthogonality relation \eqref{ortho3} yields
$$<\partial_iQ^3,R_\bot>=0.$$
We will show that
$$\underset{f\in^\bot\{Q,\partial_iQ^3\}}\inf\frac{<Lf,f>}{\|f\|_2^2}=
I>0.$$
From the proof of Lemma \ref{2vp} we have $I\geq 0$. 
Consider now a sequence of functions $f_j$, normalized in $\mathbb{L}^2$, 
which minimize $I$
$$<Lf_j,f_j>\underset{j\tend\infty}{\longrightarrow} I.$$
The gradients of $f_j$ are also bounded in $\mathbb{L}^2$, so we can 
extract a subsequence converging weakly in $\mathbb{H}^1$ to a function 
$f$
$$f_{j_n}\rightharpoonup f.$$ 
In particular,
$$<f_{j_n}^2,Q^2>\underset{n\tend\infty}{\longrightarrow}<f^2,Q^2>,$$
and it follows that $f$ is a minimizer for $I$,
$$Lf=If.$$
If $I=0$, then $f$ must be in the kernel of $L$. 
Lemma \ref{ker} ensures us that the kernel contains only the 
derivatives of $Q$, and since $f$ is orthogonal to the derivatives of 
$Q^3$, it follows that $f=0$. 
This is in contradiction with the positive $\mathbb{L}^2$ norm of $f$, 
so $I>0$.  

Therefore, since $R_\bot$ is orthogonal to $Q$ and 
to the two derivatives of $Q^3$, 
$$<LR_\bot,R_\bot>\geq I\|R_\bot\|_2^2=I(\|\tilde{R}\|_2^2-
\|R_\shortparallel\|_2^2).$$
Arguing as for the first term,
$$\|R_\shortparallel\|_2^2\leq C\|\tilde{R}\|_2^4,$$
and we finally have
$$<L \tilde{R},\tilde{R}>\geq I\|\tilde{R}\|_2^2-C\|\tilde{R}\|_2^4-
C\|\tilde{R}\|_2^3.$$
Since $\tilde{R}$ tends to $0$ in $\mathbb{L}^2$ norm, there exist a
positive constant $C$ such that for $t$ close enough to $T$,
$$<L \tilde{R},\tilde{R}>\geq C\|\tilde{R}\|_2^2.$$
For a positive number $\epsilon$,
$$<L \tilde{R},\tilde{R}>=\epsilon \left(\int|\nabla \tilde{R}|^2+
\int (1-3Q^2)\tilde{R}^2\right)+(1-\epsilon)<L \tilde{R},\tilde{R}>,$$
so, using the control of the $\mathbb{L}^2$ norm by $L$ and the
boundeness of $Q$,
$$<L \tilde{R},\tilde{R}>\geq \epsilon \left(\int|\nabla \tilde{R}|^2-
C_Q\int \tilde{R}^2\right)+(1-\epsilon)C\|\tilde{R}\|_2^2.$$
By choosing $\epsilon$ small enough to have
$$(1-\epsilon)C-\epsilon \,C_Q>0,$$
we get the existence of a positive constant $\delta$ such that 
$$<L \tilde{R},\tilde{R}>\geq \delta\|\tilde{R}\|_{\mathbb{H}^1}^2.$$
\end{proof}

Therefore the proof of Proposition \ref{dec} is complete.

\end{document}